\documentclass{article}

\usepackage{amsfonts}
\usepackage{amssymb}
\usepackage{enumerate}
\usepackage{amsmath}
\usepackage[T1]{fontenc}

\newtheorem{ttt}{Theorem}[section]
\newtheorem{llll}[ttt]{Lemma}
\newtheorem{ccc}[ttt]{Claim}
\newtheorem{eee}[ttt]{Example}
\newtheorem{fff}[ttt]{Fact}
\newtheorem{rrr}[ttt]{Remark}
\newtheorem{sss}[ttt]{Statement}
\newtheorem{ddd}[ttt]{Definition}
\newtheorem{qqq}[ttt]{Question}
\newtheorem{cccc}[ttt]{Corollary}
\newtheorem{nnn}[ttt]{Notation}

\newcommand{\bt}{\begin{ttt}}
\newcommand{\bl}{\begin{llll}}
\newcommand{\bc}{\begin{ccc}}
\newcommand{\bex}{\begin{eee}}
\newcommand{\bfa}{\begin{fff}}
\newcommand{\br}{\begin{rrr}\upshape}
\newcommand{\bs}{\begin{sss}}
\newcommand{\bd}{\begin{ddd}\upshape}
\newcommand{\bq}{\begin{qqq}}
\newcommand{\bnn}{\begin{nnn}}
\newcommand{\bcor}{\begin{cccc}}

\newcommand{\bp}{\noindent\textbf{Proof. }}

\newcommand{\et}{\end{ttt}}
\newcommand{\el}{\end{llll}}
\newcommand{\ec}{\end{ccc}}
\newcommand{\eex}{\end{eee}}
\newcommand{\efa}{\end{fff}}
\newcommand{\er}{\end{rrr}}
\newcommand{\es}{\end{sss}}
\newcommand{\ed}{\end{ddd}}
\newcommand{\eq}{\end{qqq}}
\newcommand{\ecor}{\end{cccc}}
\newcommand{\enn}{\end{nnn}}

\newcommand{\ep}{\hspace{\stretch{1}}$\square$\medskip}

\newcommand{\lab}[1]{\label{#1}}

\newcommand{\NN}{\mathbb{N}}
\newcommand{\ZZ}{\mathbb{Z}}
\newcommand{\QQ}{\mathbb{Q}}
\newcommand{\RR}{\mathbb{R}}

\newcommand{\al}{\alpha}
\newcommand{\ga}{\gamma}

\newcommand{\om}{\omega}

\newcommand{\ka}{\kappa}
\newcommand{\la}{\lambda}

\newcommand{\iA}{\mathcal{A}}
\newcommand{\iB}{\mathcal{B}}

\newcommand{\iD}{\mathcal{D}}

\newcommand{\iI}{\mathcal{I}}

\newcommand{\iM}{\mathcal{M}}

\newcommand{\iF}{\mathcal{F}}
\newcommand{\iP}{\mathcal{P}}

\newcommand{\iN}{\mathcal{N}}
\newcommand{\iT}{\mathcal{T}}

\newcommand{\RRR}{\RR^\RR}

\def\su{\subset}
\def\se{\setminus}
\def\nl{[0,1]}
\def\De{\Delta}
\def\stb{,\ldots ,}

\newcommand{\cf}{\mathrm{cf}} 

\newcommand{\solc}{{\rm sc}}
\newcommand{\non}{{\rm non}}
\newcommand{\cof}{{\rm cof}}
\newcommand{\add}{{\rm add}}

\title{A cardinal number connected to the
solvability of systems of difference equations in a given function class}

\author{M\'arton Elekes\thanks{Partially supported by Hungarian Scientific
Foundation grants no.~49786, 37758 and F 43620.} \ and Mikl\'os 
Laczkovich\thanks{Partially supported by Hungarian Scientific Foundation 
grant no.~49786.}}

\begin{document}

\maketitle 

\begin{abstract}
Let $\RR^\RR$ denote the set of real
valued functions defined on the real line.
A map $D: \RR^\RR \to \RR^\RR$ is said to be a {\it
difference operator}, if there are real numbers $a_i , b_i
\ (i=1,\ldots , n)$ such that $(Df)(x)=\sum_{i=1}^n a_i f(x+b_i)$ for
every $f\in \RR^\RR $ and $x\in \RR$. 
By a {\it system of difference equations} we mean a set of equations
$S=\{ D_i f=g_i : i\in I\}$, where $I$ is an arbitrary set of indices,
$D_i$ is a difference operator and $g_i$ is a given function for
every $i\in I$, and $f$ is the unknown function. 
One can prove
that a system $S$ is solvable
if and only if every finite subsystem of $S$ is solvable.
However, if we look for solutions belonging to a given class
of functions, then the analogous statement is no longer true. For
example, there exists a system $S$ such that every finite subsystem
of $S$ has a solution which is a trigonometric polynomial, but $S$
has no such solution; moreover, $S$ has no measurable solutions.

This phenomenon motivates the following definition. Let ${\cal F}$
be a class of functions. The {\it solvability cardinal} $\solc (\iF )$ 
of ${\cal F}$
is the smallest cardinal number $\kappa$ such that whenever $S$ is a
system of difference equations and each subsystem of $S$ of cardinality
less than $\kappa$ has a solution in ${\cal F}$, then $S$ itself
has a solution in ${\cal F}$. 
In this paper we will determine the
solvability cardinals of most function classes that occur in analysis.
As it turns out, the behaviour of $\solc (\iF )$ is rather erratic.
For example, $\solc (\text{polynomials})=3$ but $\solc (\text{trigonometric
polynomials})=\om _1$,
$\solc (\{ f: f\ \text{is continuous}\}) = \om_1$ but $\solc (\{ f: f\ 
\text{is Darboux}\}) =(2^\om )^+$, and $\solc (\RRR )=\om$.
We consistently determine the
solvability cardinals of the classes of Borel, Lebesgue and Baire measurable 
functions, and give some partial answers for the Baire class
1 and Baire class $\alpha$ functions.  
\end{abstract}

\insert\footins{\footnotesize{MSC codes: Primary 39A10, 39A70, 47B39, 26A99; 
Secondary 03E35, 03E17}}
\insert\footins{\footnotesize{Key Words: systems of difference
equations, cardinal invariants}}

\section{Preliminaries}

Difference operators occur in various branches of analysis. For
example, it is shown in \cite{La2} that the existence of certain
types of liftings is closely related to the solvability of systems of
difference equations. Among others, 
it is obtained from results on the solvability of infinite
systems of difference equations that there exists a linear
operator from the bounded real functions into the set of measurable
real functions that fixes the bounded measurable functions and
commutes with any prescribed countable set of translations
\cite[Theorem 3.3]{La2}. On the other hand, there is no such linear
operator from the space of all complex valued functions defined on $\RR$
into the space $L^0$ of measurable functions; see 
\cite[Theorem 5.1 and 5.2 ]{La2}.

The goal of this paper is to give necessary
conditions under which systems of difference equations have solutions
belonging to a given function class.

\bnn
{\rm 
Let $\RR^\RR$ denote the set of real
valued functions defined on the real line. The classes of polynomials and
trigonometric polynomials are denoted by $\iP$ and $\iT\iP$.
For every set $H$ we shall denote by $\chi _H$ and $|H|$ the
characteristic function and the cardinality of $H.$ We denote the symmetric
difference of the sets $A$ and $B$ by $A\Delta B$. 
If $A,B\su \RR$ and $x\in \RR$ then we shall write
$A+B=\{ a+b:a\in A,\ b\in B\}$ and $A+x=\{ a+x:a\in A\} .$
If $A\subset \RR$
then $\langle A\rangle$ denotes the additive group generated by $A.$
The symbols $\kappa^+$ and $\cf(\ka)$ denote the successor cardinal and the
cofinality of the cardinal $\kappa$.}
\enn

\bd 
A \emph{difference operator} is a mapping $D:\RR^\RR\to\RR^\RR$ of the
form
\[
(Df)(x)=\sum_{i=1}^n a_i f(x+b_i),
\]
where $a_i$ and $b_i$ are real numbers. The set of difference
operators is denoted by $\iD$.
\ed

\bd
For $b\in\RR$ the difference operators $T_b$ and $\Delta_b$ are defined by
\[
(T_b f)(x) = f(x+b)\ \ \ (x\in\RR), \ \ \text{and}
\]
\[
(\Delta_b f)(x) = f(x+b) - f(x) \ \ \ (x\in\RR).
\]
\ed

\bd
A \emph{difference equation} is a functional equation
\[
Df=g,
\]
where $D$ is a difference operator,
$g$ is a given function and $f$ is the unknown. A \emph{system of
difference equations} is 
\[
D_i f=g_i \ \ (i\in I),
\]
where $I$ is an arbitrary set of indices. More formally, by
a system of difference equations we mean a set
$S\su \iD \times \RRR .$
A function $f:\RR\to\RR$ is a solution
to $S$ if $Df= g$ for every $(D,g)\in S.$
\ed

It was proved in \cite[Thm.~2.2]{La2} that a system of difference
equations is solvable iff each
of its finite subsystems is solvable. However, if we are interested
in solutions belonging to a given subclass of $\RR ^\RR$ then this result is no longer
true. This motivates the following. 

\bd
Let $\mathcal{F}\subset \RR^\RR$ be a class of real functions. The
\emph{solvability cardinal of $\mathcal{F}$} is the minimal cardinal
$\solc(\mathcal{F})$ with the property that if every subsystem of size
less than $\solc(\mathcal{F})$ of a system of difference equations has a
solution in $\mathcal{F}$, then the whole system has a solution in
$\mathcal{F}$.
\ed

For example, $\solc(\RR^\RR)\leq\omega$ is a reformulation of the above
cited result.
The next statement shows that the cardinal $\solc(\iF)$ actually exists,
and also provides an upper bound.

\bfa\lab{fact}
For every $\iF\subset\RRR$ we have $\solc(\iF)\leq (2^\omega)^+$.
\efa

\bp
Note that the cardinality of $\iD$ is $2^\omega$. Suppose
$\iF\subset\RRR$, $S$ is system of difference equations, and 
every subsystem of $S$ of cardinality at most $2^\omega$ is
solvable in $\iF$. In particular, every pair of equations of $S$ is
solvable, hence for every $D\in\iD$ there is at most one $g\in\RRR$
such that $(D,g)\in S$. Therefore the cardinality of $S$ is at most
$2^\omega$, and we are done.
\ep

We may add the following trivial estimate.

\bfa\lab{fact2}
For every $\iF\subset\RRR$ we have $\solc(\iF)\leq |\iF |^+$.
\efa

\bp
Let $S$ be a system of difference equations such that
every subsystem of $S$ of cardinality at most $|\iF |$ is
solvable in $\iF$. Suppose $S$ is not solvable in $\iF$.
Then for every $f\in \iF$ there is a $(D_f ,g_f )\in S$ such that $D_f f\ne
g_f .$ Then $S'=\{ (D_f ,g_f ): f\in \iF \}$ has no solution in $\iF$
and $|S'|\le |\iF |,$ a contradiction.
\ep

\br
Fact \ref{fact2} can be improved if we take into consideration
the product topology on $\RRR .$ Namely, if $Df\ne g$ for some
$f\in \RRR$ then $f$ has a neighbourhood $U$ in the product topology
such that $Df'\ne g$ for every $f'\in U.$ Combining this observation with
the proof of Fact \ref{fact2} we obtain the estimate $\solc (\iF )\le 
L(\iF )^+ ,$ where $L(X)$ is the Lindel\"of number of the topological space $X$;
that is, the smallest cardinal $\ka$ such that each open cover of $X$ 
contains a subcover of cardinality at most $\ka .$ This sharper
inequality implies Fact \ref{fact} since the space $\RRR$ has a base
of cardinality $2^\om$ and thus $L(X)\le 2^\om$ for every subspace $X\su \RRR
.$
\er

It is natural to ask whether or not every cardinal $2\le \ka \le (2^\om )^+$
equals $\solc (\iF )$ for some $\iF \su \RRR .$ As we shall see
in Theorem \ref{t:a}, $\om$ is such a cardinal. 
The following
result gives a positive answer for successor cardinals.

\bt \label{t:succ}
For every cardinal $1\le \ka \le 2^\om$ there exists an $\iF \su \RRR$
such that $\solc (\iF )=\ka ^+ .$
\et

\bp
Let $B\su \RR$ be linearly independent over the rationals with $|B|=\ka$.
For every $b\in B$ we denote by $f_b$ the characteristic function of the
group $\langle B\se \{ b\} \rangle .$ Then $f_b$ is periodic mod 
each element of $B\se \{ b\} $,
but $f_b$ is not periodic mod $b.$

We claim that the solvability cardinal of the class $\iF =\{ f_b :b\in B\}$
equals $\ka ^+ .$ The inequality $\solc (\iF )\le \ka ^+$ is clear from
Fact \ref{fact2}. In order to prove $\solc (\iF )\ge \ka ^+$ we have to 
construct a system $S$ such that every subsystem $S'\su S$ of size less
than $\ka$ is solvable in $\iF ,$ while $S$ is not.
We show that $S=\{ (\De _b ,0):b\in B\}$ is such a system.
If $b\in B$ then, as $f_b$ is not periodic mod $b,$ we have $\De _b f_b \ne 0$ 
showing that $S$ is not solvable in $\iF .$ On the other hand, if $S'$ is a 
proper subsystem of $S$ and $(\De _b ,0) \notin S'$ then $f_b$ solves
$S'$ completing the proof.
\ep

\bq
Is it true (in ZFC) that for every $2\le\kappa\le (2^\om)^+$ there exists an
$\iF \su \RRR$ such that $\solc(\iF) = \ka$? Is there (in ZFC) an
$\iF$ with $\solc(\iF) = 2^\om$? Is it consistent with ZFC that
$\solc(\iF)$ can be an uncountable limit cardinal?
\eq
 
In the first part of the paper (Sections \ref{s:a}, \ref{s:b},
\ref{s:D}, and \ref{s:s})
we determine the exact value of $\solc(\iF)$ for several classes $\iF$; see
Theorems \ref{t:a}, \ref{t1}, \ref{t:b}, \ref{D},
Corollaries \ref{c:1}, \ref{c:2} and
Theorem \ref{t:p}. 
As it turns out, the behaviour of $\solc (\iF )$ is rather erratic.
For example, $\solc (\iP )=3,$ but $\solc (\iT \iP )=\om _1 ;$
$\solc (\RRR )=\om ,$ but $\solc (\{ f: f\ \text{is Darboux}\} =(2^\om )^+ .$

In the second part (Sections \ref{s:Borel}, \ref{s:Baire}
and \ref{s:LB}) we give estimates of 
$\solc(\iF)$ for subclasses of Borel, Lebesgue and Baire measurable functions.
The estimates for Borel, Lebesgue and Baire measurable functions
provide the exact value of the solvability cardinals
consistently. The result $\solc(\{f: f\textrm{ is
Lebesgue measurable}\})$ $> \omega _1$ answers Problem 3 of \cite{La2}.

\section{Arbitrary functions}\lab{s:a}

The nontrivial direction of the next theorem was proved in
\cite[Thm.~2.2]{La2}, but we reformulate this result using the
notation introduced in the present paper.

\bt\lab{t:a}
$\solc(\RRR)=\omega$.
\et

\bp
$\solc(\RRR)\leq\omega$ is \cite[Thm.~2.2]{La2}, so we only need to show
that $\solc(\RRR)\neq n$ for every $n\in\NN$. Let $n\geq 2$,
let $a_1, \ldots, a_{n-1} \in \RR$ be 
linearly independent over the rationals, and put $a_n
=-\sum_{i=1}^{n-1} a_i$. Then any $n-1$ of the numbers $a_1 ,\ldots
,a_n$ are linearly independent over the rationals.
Define the following system of $n$ equations:
\[
\Delta_{a_i} f = 1, \ \ i=1,\ldots,n.
\]
It is easy to see that each subsystem of cardinality at most $n-1$ is
solvable (consider the factor group of $\RR$ modulo the additive
group generated by the corresponding 
linearly independent $a_i$'s). On the other hand, if $f$ were a
solution to the whole system, then $f(0)+n=f(a_1+\ldots+a_n)=f(0)$
would hold, which is impossible. This shows
$\solc(\RRR)>n$ and, as $n$ was arbitrary, the
proof is complete.
\ep

\section{Bounded functions}\lab{s:b}

It is well known that the difference operators form an algebra
under the operations $(A+B)f=Af +Bf, \ (c\cdot A)f =c\cdot Af$ and
$(AB)f=A(Bf).$

\bd
We say that the difference equation $(D,g)$ is \emph{deducible from} 
the system $S$ if there are $A_1,\ldots, A_n \in \iD$ and $(D_1,g_1),
\ldots,(D_n,g_n) \in S$ such that $(D,g) = (\sum_{i=1}^n A_i D_i,
\sum_{i=1}^n A_i g_i)$.
\ed

\bt \label{t1}
Let $K>0$ be a real number. Then $\solc(\{f \in \RRR: |f|\leq K\}) = \omega$.
\et

\bp 
We may assume $K=1$. First we show $\solc(\{f \in \RRR: |f|\leq 1\}) \leq
\omega$. The proof is a modification of the proof of
\cite[Thm.~2.1]{La2}, the new ingredient is the Hahn-Banach
Theorem. Let $S$ be a system such that all finite subsystems are
solvable by functions of absolute value at most $1$. Define
\[
\iA=\{D\in\iD : \exists g \ (D,g) \textrm{ is deducible from } S
\}.
\]
Then $\iA$ is a linear subspace of $\iD$. Put
\[
L(D) = g(0) \ \ \ (D\in\iA),
\]
where $(D,g)$ is deducible from $S$. 
Clearly, if $(D,g)$ is
deducible from $S$ then it is also deducible from a finite subsystem
of $S$, hence it is solvable. Moreover, any pair of equations
deducible from $S$ has a common solution. Therefore the map $L:\iA \to
\RR$ is well defined.
Note that $L$ is clearly linear.

Now we define a norm on $\iD$. 
It is easy to see that every
$D\in\iD$ has a unique representation of the form $D = \sum_{i=1}^n a_i
T_{b_i}$, where the $a_i$'s are nonzero and the $b_i$'s are
different. Using this representation set
\[
||D|| = \sum_{i=1}^n |a_i|.
\]
The function $||.||:\iD\to\RR$ is easily seen to be a norm. 

We claim that for every $D\in\iA$ we have $|L(D)| \leq ||D||$. Let
$(D,g)$ be deducible from $S.$ Then there is a function $f$ such that
$|f|\le 1$ and $Df=g.$ If $D=\sum_{i=1}^n a_i T_{b_i}$ then
$$|L(D)|=|g(0)|= \left| \sum_{i=1}^n a_i f(b_i ) \right| \leq
\sum_{i=1}^n |a_i| \cdot 1 =||D||.$$
Hence by the Hahn-Banach Theorem (see e.g.~\cite[Thm. 3.3]{Ru1}) there
exists a linear map
$L^*:\iD\to\RR$ extending $L$ such that
\[
\left|L^*(D)\right| \leq ||D|| \textrm{ for every } D\in\iD.
\]
We claim that the function defined by
\[
f(x) = L^*(T_x) \ \ \ (x\in\RR)
\]
is a solution to $S$ such that $|f|\leq 1$. This last inequality is
obvious, as $|f(x)| = \left|L^*(T_x)\right| \leq ||T_x|| =1$. So we
need to prove that $f$ solves $S$. First we show that
\begin{equation}\lab{param}
(Df)(0)=L^*(D) \textrm{ holds for every } D \in \iD.
\end{equation}
Since $L^*$ is linear, it is enough to check this for $D=T_x \ (x\in
\RR)$. 
Now $(T_x f) (0)= f(x)=L^*( T_x ) $ by the definition of $f$, which
proves (\ref{param}). Let $(D,g)\in S$ and $x\in \RR$ be given.
Then $T_x D\in \iA ,$ and thus (\ref{param}) and the definition of $L$
imply
\[
(Df)(x)=( T_x D f) (0)= L^* ( T_x D ) = L ( T_x D ) = ( T_x
g ) (0)=g(x).
\]

Now we prove $\solc(\{f \in \RRR: |f|\leq 1\}) \geq \omega$. Let $n\geq
2$ be an integer and let $a_1,\ldots,a_n$ be linearly independent
reals. Define a system as follows.
\[
\Delta_{a_i} f =
\frac{2}{n-1}\ \chi_{ \langle \{a_1,\ldots,a_{i-1},a_{i+1},a_n\}
\rangle }, \ \
(i=1,\ldots,n).
\]
A simple induction shows that if $f$ solves the whole system, then 
$f(a_1+\ldots+a_n) - f(0) = \frac{2n}{n-1} >2$, hence $|f|\leq 1$ cannot
hold.

On the other hand, let $J\subset \{1,\ldots,n\}$ be a set of at most $n-1$
elements. Every $x\in\langle\{a_1,\ldots,a_n\}\rangle$ can be uniquely
written in the form $x = k_1(x)a_1+\ldots+k_n(x)a_n$, where the
$k_i(x)$'s are integers. 
Define
\[
f(x) = \begin{cases}  0 & \text{if $x\notin\langle\{a_1,\ldots,a_n\}\rangle$} \cr
               -1+\frac{2}{n-1}\ |\{i\in J : k_i(x) > 0 \}| & \textrm{if
               $x\in\langle\{a_1,\ldots,a_n\}\rangle$.}
\end{cases}
\]
Clearly, $|f|\leq 1$. It is easy to see that $f$ solves the $i^{th}$ 
equation for every $i\in J$, which yields $\solc(\{f \in \RRR: |f|\leq
1\}) > n$. As $n$ was arbitrary, the proof is complete.
\ep

In contrast to Theorem \ref{t1} we have the following.

\bt\lab{t:b}
$\solc(\{f \in \RRR:\  $f$ \textrm{ is bounded}\}) = \omega_1$.
\et

\bp
First we prove $\solc(\{f \in \RRR:\  f \textrm{ is bounded}\}) \leq
\omega_1$. Let $S$ be a system such that every countable subsystem of
$S$ is solvable by a bounded function. For a countable $S'\subset S$
let $K_{S'}$ be the minimal integer for which $S'$ has a solution in
$\{f \in \RRR: |f|\leq K_{S'}\}$. The set $\{K_{S'} : S'\subset S, \ 
|S'|\leq \omega \}$ is bounded in $\NN$, otherwise we could easily
find a countable subsystem of $S$ with no bounded solutions. Fix an
upper bound $K$ of the above set. Then every countable, in
particular, every finite subsystem of $S$ is solvable in $\{f \in
\RRR: |f|\leq K\}$, hence by the previous theorem  $S$ is solvable in
$\{f \in \RRR: |f|\leq K\}$, hence $S$ has a bounded solution.

Now we prove $\solc(\{f \in \RRR:\  f \textrm{ is bounded}\}) >
\omega$. Similarly to the previous theorem, let $a_1,a_2,\ldots$ be a
linearly independent sequence of reals. Define a system by
\[
\Delta_{a_i} f =
\chi_{\langle \{a_1,\ldots,a_{i-1},a_{i+1},\ldots\} \rangle}, \ \
(i\in\NN^+).
\]
A simple induction shows that if $f$ solves the whole system, then 
$f(a_1+\ldots+a_n) - f(0) = n$ for every $n$, hence $f$ cannot be bounded.

On the other hand, let $J\subset \NN^+$ be a finite set. Every
$x\in\langle\{a_1,a_2,\ldots \}\rangle$ can be uniquely
written in the form $x = k_1(x)a_1+k_2(x)a_2+\ldots$, where the
$k_i(x)$'s are integers, and only finitely many of them are nonzero. 
Similarly to the proof of the previous theorem one can check that
\[
f(x) = \begin{cases} 0 & \text{if $x\notin\langle\{a_1,a_2,\ldots\}\rangle$} \cr
               |\{i\in J : k_i(x) > 0 \}| & \text{if
               $x\in\langle\{a_1,a_2,\ldots\}\rangle$ }
\end{cases}
\]
is a bounded solution to the finite subset of $S$ corresponding to $J$.
\ep

\section{Darboux functions}\lab{s:D}

\bt\lab{D}
$\solc(\{f : f \textrm{ is Darboux}\}) = (2^\omega)^+$.
\et

\bp
$\solc(\{f : f \textrm{ is Darboux}\}) \le (2^\omega)^+$ follows from
Fact \ref{fact}. In order to prove the other inequality we have to
construct a system $S$ such that every subsystem of cardinality less
than continuum is solvable by a Darboux function but $S$ has no
Darboux solution. We define $S$ as
\[
\Delta_b f = \Delta_b \chi_{\{0\}} \ \ \ (b\in\RR).
\]
The whole system clearly has no Darboux solution, for if $f$ is a
solution to $S$ then there exists a $c\in\RR$ such that $f =
\chi_{\{0\}} +c$, which is not Darboux. On the other hand, let $S'$ be
a subset of $S$ such that $|S'| <2^\om$, and let $B\su\RR$ be
the corresponding set of indices with $|B| <2^\om$. 
By enlarging $B$ if necessary, we may assume that $B$ is an
additive subgroup of $\RR$, and also that $B$ is dense. 

As $|B| <2^\om$, the factor group $\RR / B$ consists of $2^\om$
cosets. Fix a bijection $\varphi : \RR/B \to \RR$ and define
\[
f(x) = \varphi(B+x) + \chi_{\{0\}}(x) \ \ \ (x\in\RR).
\]
As $B$ is dense, $f$ attains every value 
on every interval, hence it is Darboux. In addition, it is easy to see
that $f$ solves $S'$.
\ep

\br
The same system can be used to demonstrate that for the class $\iF$ of
functions with connected graphs we also have $\solc(\iF) =
(2^\omega)^+$. With a more elaborate version of the argument above it
can be shown that if $|B| < 2^\om$ then the
system $\{(\Delta_b, \Delta_b \chi_{\{0\}}) : b\in B \}$ has a
solution with a connected graph.
\er

\section{Subclasses of Lebesgue measurable and Baire measurable
functions}\lab{s:s}

In this section our aim is to prove that $\solc (\iF )=\om_1$
for many classes including the classes of 
trigonometric polynomials, 
continuous functions, Lipschitz functions, $C^n$, $C^\infty$, analytic
functions, derivatives, approximately continuous functions etc.

Let $\iN$ denote the $\sigma$-ideal of Lebesgue nullsets of $\RR$ and
$\iM$ denote the $\sigma$-ideal of first category (=
meager) subsets of $\RR$. 
In the rest of the section let $\iI$ stand for either $\iN$
or $\iM$. The term $\iI$-almost everywhere will be abbreviated by $\iI$-a.e.
Instead of 'Lebesgue measurable' and 'with the
Baire property' we will use the term \emph{$\iB_{\iI}$-measurable},
where $\iB_{\iI}$ is the $\sigma$-algebra generated by the Borel sets
and $\iI$.
 
First we show that if we do not distinguish between $\iI$-almost
everywhere equal functions, then
the value of this modified solvability cardinal is at most $\omega_1$
for all subclasses of both
Lebesgue measurable functions and functions with the property of
Baire. 

\bt\lab{ae}
Let $\iF \su \iB _{\iI},$ and
suppose that for every countable subsystem $S'$ of a system of
difference equations $S$ there exists an
$f' \in \iF $ such that $Df'=g \ \iI$-a.e.~for every $(D,g) \in
S'$. Then there is an $f\in \iF$ such that
$Df=g \ \iI$-a.e.~for every $(D,g) \in S$.
\et

\bp
Let $S$ be a system satisfying the assumptions.
Every $D\in\iD$ can be written in a unique way as
$D = \sum_{i=1}^n a_i T_{b_i}$. Define $\varphi : \iD \to
\bigcup_{n\in\NN} \RR^{2n}$ by
\[
\varphi(D) = (a_1,\ldots,a_n,b_1,\ldots,b_n) \ \ (D\in\iD).
\]
Set $S_n = \{(D,g)\in S : D \textrm{ has $n$ terms}\}$. For every
$n\in \NN$ choose a countable $S_n' \subset S_n$ such that $\{\varphi(D)
: (D,g)\in S_n'\}\subset \RR^{2n}$ is dense in $\{\varphi(D)
: (D,g)\in S_n\}\subset \RR^{2n}$. 

Let $f\in \iF$ be a
function '$\iI$-a.e.' solving $\bigcup_{n\in\NN} S_n'$. We claim
that it '$\iI$-a.e.' solves the whole $S$. Let $(D,g) \in S_n$, and 
choose $(D_i,g_i) \in S_n'$ such that $\varphi(D_i) \to \varphi(D)$ in
$\RR^{2n}$. 

Suppose first $\iI = \iN$. It is well
known that for every measurable $h$ if
$t_n \to 0$ 
$(n\to\infty)$ then $T_{t_n} h \to h$ in measure (which means 
that it converges in
measure on every bounded interval; see e.g.~\cite{Ru2} or
\cite{Ha} for the definitions and basic facts). Hence
$D_i h \to D h$ in measure. Let $f'$ be an a.e.~solution to
$\bigcup_{n\in\NN} S_n' \cup\{(D,g)\}$. Then
$$Df=\lim_{i\to \infty} D_i f=\lim_{i\to \infty} g_i =\lim_{i\to
\infty} D_i f' =Df' =g$$
a.e., where $\lim$ stands for limit in measure. 


Suppose now $\iI = \iM$. We claim that for every $h$ with the Baire 
property if $t_n \to 0$ $(n\to\infty)$ then $T_{t_n} h \to h$ 
pointwise on a residual set. Indeed, if $H$ is a
residual set on which $h$ is continuous then $H \cap \bigcap_{n\in\NN}
(H-t_n)$ is such a set. Therefore
$D_i h \to D h $ pointwise on a residual set.
Let $f'$ be an $\iM$-a.e.~solution to
$\bigcup_{n\in\NN} S_n' \cup\{(D,g)\}$. Then 
$$Df=\lim_{i\to \infty} D_i f=\lim_{i\to \infty} g_i =\lim_{i\to
\infty} D_i f' =Df' =g$$
on a residual set. \ep

\bt \label{t2}
Let $\iF \su \tilde{\iF} \su \iB _{\iI},$ where $\tilde{\iF}$ is 
a translation invariant linear subspace of $\iB _{\iI}$
such that whenever $f\in \tilde{\iF}$ and $f=0 \ \iI$-a.e.~then 
$f=0$ everywhere. Then $\solc (\iF )\le \omega _1 .$
\et

\bp
Suppose that every countable subsystem of
$S$ has a solution in $\iF .$ Then obviously $g\in \tilde{\iF}$
whenever $(D,g)\in S.$ 

By Theorem \ref{ae},
there is an $f\in \iF$ such that $Df=g$ $\iI$-a.e.~for every $(D,g)\in
S.$ Since $Df-g\in \tilde{\iF}$ and $Df-g=0$ $\iI$-a.e., we have $Df=g,$ 
which proves $\solc (\iF )\le \omega _1 .$ \ep

It is clear that the class $C(\RR )$ of continuous functions satisfies
the conditions imposed on $\tilde{\iF} .$ The same is true for the
classes of derivatives and approximately continuous functions (see
\cite{Br}). 

We shall denote by $\iT\iP$ the set of trigonometric polynomials.

\bt \label{t4}
If $\iT\iP \su \iF \su \iB _{\iI}$ then
$\solc (\iF )\ge  \omega _1 .$
\et

\bp
We shall construct a system $S$ such that every
finite subsystem of $S$ has a solution which is a trigonometric polynomial,
but $S$ itself does not have a $\iB _{\iI}$-measurable solution.
We shall repeat the construction of
\cite[Thm.~4.4]{La2} with a small modification.

Let $C(x)=\cos 2\pi x $ and
$E_{j,n} (x)=\De_{2^{-n}} C\left( 2^j x\right),$ 
then $E_{j,n} \in \iT \iP$ for every $j,n\in \NN .$
Also,
$E_{j,n}=0$
if $j\ge n$ and, if $j<n$ then $E_{j,n}$ is a continuous function
periodic mod 1 with finitely many roots in $[0,1].$

Let $c_j \ (j=0,1,\ldots )$ be a sequence
of real numbers, and consider the system $S$ of the equations
$$\De_{2^{-n}} f=h_n ,\ \
{\rm where}\ \	h_n = \sum_{j=0}^{n-1} c_j E_{j,n}  \ \
(n=1,2, \ldots ).$$
Then the trigonometric polynomial
$\sum_{j=0}^{n-1} c_j  C\left( 2^j x\right)$ is a solution to the first
$n$ equations of $S.$ On the other hand, we shall choose the numbers $c_j$
in such a way that $S$ does not have $\iB _{\iI}$-measurable solutions.

First suppose $\iI =\iN .$
If $f:\RR \to \RR$ is measurable then the sequence
of functions $\De_{2^{-n}} f$ converges to zero in measure on $\nl .$
Therefore, if $S$ has a measurable solution, then $h_n$ should
converge to zero in measure on $[0,1] .$ But we can prevent this by a suitable
choice of the sequence $c_j .$ We shall define $c_j$ inductively. If $c_j$
has been defined for every $j<n-1,$ then we choose $c_{n-1}$ so large
that $\la (\{ x\in [0,1] : |h_n (x)|>1\} )>1/2$ holds. This is possible,
since $E_{n-1,n} \ne 0$ a.e.~in $[0,1].$
Therefore, with this
choice, $h_n$ does not converge (in measure) to zero on $\nl ,$ and
thus $S$ cannot have measurable solutions. 

Next suppose $\iI =\iM .$
If $f:\RR \to \RR$ is Baire measurable then the sequence
of functions $\De_{2^{-n}} f$ converges to zero pointwise on a residual subset
of $\nl .$ Again, we shall choose the constants $c_j$ such
that $h_n \not\to 0$ on a second category set.
Namely, we shall define $c_j$ in such a way that each function $h_n$
satisfies
the following condition: for every interval $I \su \nl$ of length
$1/n$ the inequality $|h_n |>1$ holds on a subinterval of $I.$ (In the
course of the proof 
by an interval we shall mean a closed nondegenerate interval, and
by $|I|$ we shall mean the
length of the interval $I$.)

We put $c_0 =1.$ Then $h_1(x) = C(x+\frac{1}{2})-C(x) = -2\cos2\pi x$ has the
required property with $n=1,$ since there is a subinterval of $\nl$ on which
$|h_1|>1$. 
Let $n>1$ and suppose that $c_0 \stb c_{n-2}$ have been chosen. Since
$E_{n-1,n}$ only has a finite number of roots in $\nl ,$ the function
$$h_n =\left( \sum_{j=0}^{n-2} c_j E_{j,n} \right) +c_{n-1} E_{n-1,n}$$
clearly has the required property if $c_{n-1}$ is large enough.

We show that the set $A=\{ x\in \nl : h_n (x)\to 0\}$ is not residual.
Suppose the contrary, and let $\bigcap _{k=1}^\infty U_k \su A,$
where each $U_k$ is dense open. Let $I_1 \su U_1$ be an
interval. If $1/n_1 < |I_1 |$ then there is a subinterval $J_1 \su I_1$
such that $|h_{n_1}| >1$ on $J_1 .$ Since $U_2$ is dense open, there
is an interval $I_2 \su U_2 \cap J_1 .$ If $1/n_2 < |I_2 |$ then there is a subinterval 
$J_2 \su I_2$
such that $|h_{n_2}| >1$ on $J_2 .$ Continuing this process we find the nested 
sequence of intervals $J_k$ such that $\bigcap_{k=1}^\infty J_k \su 
\bigcap_{k=1}^\infty 
U_k \su A$ and $|h_{n_k} |>1$ on $J_k$ for every $k.$ This implies
$h_n (x) \not\to 0$ for every $x\in \bigcap_{k=1}^\infty J_k $
which contradicts $x\in A.$ \ep

\bcor\lab{c:1}
Suppose
$\iT \iP \su \iF \su \tilde{\iF} \su \iB _{\iI},$ where $\tilde{\iF}$ is 
a translation invariant linear subspace of $\iB _{\iI}$
such that whenever $f\in \tilde{\iF}$ and $f=0 \ \iI$-a.e.~then $f=0$
everywhere. Then $\solc (\iF )= \omega _1 .$
\ecor

It is clear that the class $C(\RR )$ of continuous functions satisfies
the conditions imposed on $\tilde{\iF} .$ The same is true for the
classes of derivatives and approximately continuous functions (see
\cite{Br}). 
Thus we have the following.

\bcor\lab{c:2}
If $\iF$ equals any of the classes
$\iT\iP,\ C(\RR),$ the class of Lipschitz functions, 
$C^n(\RR), \ C^\infty(\RR)$, the class of
real analytic functions, 
derivatives, approximately continuous
funcions, then $\solc (\iF )=\omega _1 .$
The same is true for the subclasses $\{ f\in \iF : f\ \textrm{is 
\ bounded}\}$ where $\iF$ is any of the classes listed above.
\ecor

We remark that the class $\iP$ of polynomials behaves quite differently
from $\iT \iP .$ Indeed,
\cite[Thm.~4.5]{La2} states that $\solc(\iP) \leq 3$. 
Since $\solc(\iP) \geq 3$ is obvious, we have the following.

\bt\lab{t:p}
$\solc(\iP) = 3$.
\et

\section{Borel functions}\lab{s:Borel}

First we prove an auxiliary lemma.

\bl\lab{l:perfects}
There exist non-empty perfect subsets $\{P_\al : \al<2^\om\}$ of $\RR$ and
distinct real numbers $\{p_\al : \al<2^\om\}$ such that
\[
(P_\al+G_{\al+1}) \cap (P_\beta+G_{\beta+1}) = \emptyset \ \ \ (\al\neq\beta),
\]
and for every $\al<2^\om$
\[
(P_\al+g_1) \cap (P_\al+g_2) = \emptyset \ \ \ (g_1, g_2\in G_{\al+1}, \
g_1\neq g_2),
\]
where $G_\al =\langle \{p_\beta :
\beta<\al\} \rangle $.
\el

\bp
Let $P\subset\RR$ be a non-empty perfect set that is linearly
independent over the 
rationals (see e.g.~\cite{vN} or \cite{My}). We can choose
nonempty perfect sets $P_\al\subset P$ and $p_\al\in P$ 
($\al<2^\om$) such that $P_\al \cap P_\beta = \emptyset$ for every
$\al\neq\beta$ and such that $p_\al\notin P_\beta$ for every $\al, \beta <
2^\om$. It is a straightforward calculation to check that all the
requirements are fulfilled.
\ep

\bt\lab{Borel}
$\solc(\{f : f \textrm{ is Borel}\}) \geq \omega_2$.
\et

\bp
Let $P_\al$ and $p_\al$ be as in the previous lemma.
For every $\al<\om_1$ let $B_\al\subset P_\al$ be a Borel set of class $\al$
(that is, not of any smaller class). Define $A_\al = B_\al + G_\al ,$ and 
consider the system
of difference equations:
\[
\Delta_{p_\al} f = \Delta_{p_\al} \left(\sum_{\beta<\om_1}
\chi_{A_\beta}\right) \ \ \  (\al<\om_1). 
\]
Note that the $A_\beta$'s are disjoint.
We claim that every countable subsystem of this
system has Borel solution, but the whole system does not.

To prove the first statement we have to check that for every $\al<\om_1$ the
first $\al$ equations have a common Borel solution. We show that the Borel
function 
\[
\sum_{\beta\leq\al} \chi_{A_\beta}
\]
will do.
If $\ga<\beta$ then $A_\beta$ is periodic mod $p_\ga$, so
$\Delta_{p_\ga}\chi_{A_\beta}=0$. Therefore, in view of the properties
required in Lemma \ref{l:perfects}, we obtain that for $\ga<\al$
\[
\Delta_{p_\ga} \left(\sum_{\beta<\om_1}\chi_{A_\beta}\right) =
\Delta_{p_\ga} \left(\sum_{\beta\leq\al}\chi_{A_\beta}\right),
\]
which proves this part of the claim.

In order to show that the whole system has no Borel solution it is
sufficient to check that the functions on the right hand side of the
equations are of unbounded Baire class. But this is not hard to see, as
$\Delta_{p_\al} (\sum_{\beta<\om_1}\chi_{A_\beta})$ restricted to $P_\al$
equals $-\chi_{B_\al}$.
\ep 

Using Fact \ref{fact} we obtain the following.

\bcor\lab{borelomketto}
The Continuum Hypothesis implies that $\solc(\{f : f \textrm{ is
Borel}\}) = \om_2 = (2^\om)^+$.
\ecor

\bq
Can we omit the use of the Continuum Hypothesis? Is it true
that $\solc(\{f : f \textrm{ is Borel}\}) = \om_2$? Is it true
that $\solc(\{f : f \textrm{ is Borel}\}) = (2^\om)^+$?
\eq

\br\lab{rem}
In order to prove $\solc(\{f : f
\textrm{ is Borel}\}) = \om_2$ it would be sufficient to prove $\solc(\{f
: f \textrm{ is Baire class}~\al \}) \leq \om_2$ for every $\al<\om_1$.
Indeed, assume that every
subsystem of cardinality at most $\om_1$ of a system has a
Borel solution. Let us assign to every such subsystem the minimal
$\al<\om_1$ for which it has a Baire class~$\al$ solution. We claim
that the set of these $\al$'s is bounded in $\om_1$. Otherwise,
the union of
$\om_1$-many appropriate subsystems would itself be a subsystem of cardinality
$\om_1$ without a Borel solution, which proves our statement.

So if every subsystem of cardinality at most $\om_1$ of a system has a Borel
solution, then there exists an $\al<\om_1$ such that every such subsystem
has a Baire class~$\al$ solution. 
\er

\br
For $2\leq\al<\om_1$ the idea of the proof of Theorem \ref{Borel} probably
gives $\solc(\{f : f \textrm{ is Baire class}~\al \}) \geq
\om_2$. If we had an appropriate notion of rank for Baire
class~$\al$ functions, sharing the properties of the well known ranks
on Baire class~1, it would yield
$\solc(\{f : f \textrm{ is Baire class}~\al \}) \geq \om_2$.
Unfortunately, according to \cite{Kec2} no such rank is known.

For Baire class~1 these ranks exist, but do not give
$\solc(\{f : f \textrm{ is
Baire class}~1 \})$ $\geq \om_2$. The proof breaks down, as
$\sum_{\beta\leq\al}\chi_{A_\beta}$ is not Baire class~1. 
\er

\bq
Is there a rank on Baire class~$\alpha$ with the usual properties?
\eq

Remark \ref{rem} shows why we are particularly interested in the
solvability cardinals of the individual
Baire~$\al$ classes. The simplest case, namely $C(\RR)$ is solved
already. So we take one step further.

\section{Baire class~1 functions}\lab{s:Baire}

It is clear from Theorem \ref{t4} that $\solc(\{f : f
\textrm{ is Baire class}~1 \}) \ge  \om_1 $.
As opposed to the case $2\leq\al<\om_1$ we conjecture that, in fact,
$\solc(\{f : f \textrm{ is Baire class}~1 \}) = \om_1$. 
Unfortunately, we only can prove this in
a special case. What makes this case interesting is that it covers the
usual situation in which every difference operator $D$ is of the form
$D=\Delta_b$.

First we need two lemmas.

\bl
Let $a,b\in\RR\setminus\{0\}$. The solutions to the equation
\[
f(x+b)-af(x)=0
\]
are the functions of the form
\[
f(x)=\varphi(x)(|a|^{1/b})^x,
\]
where $\varphi$ is an arbitrary function periodic mod $b$ if $a>0$, and an
arbitrary function anti-periodic mod $b$ (that is,
$\varphi(x+b)=-\varphi(x)$ for every $x\in\RR$) if $a<0.$

In addition, $f$ is Baire class~1 iff $\varphi$ is Baire class~1.
\el

\bp
Straightforward calculations.
\ep

\bl\lab{l:ab}
Let $a_1,a_2,b_1, b_2\in\RR\setminus \{0\}$. Suppose that the equations
$f(x+b_1)-a_1f(x)=0$ and $f(x+b_2)-a_2f(x)=0$ have a common Baire class~1
solution which is not identically zero. Then $|a_1|^{1/b_1} =
|a_2|^{1/b_2}$.
\el

\bp
Suppose this is not true. Then by the previous lemma there exist two
Baire class~1 functions $\varphi_1$ and $\varphi_2$ such that
$$
\varphi_1(x)(|a_1|^{1/b_1})^x = \varphi_2(x)(|a_2|^{1/b_2})^x,
$$
where
$\varphi_1$ and $\varphi_2$ are periodic (or anti-periodic) mod $b_1$ and
$b_2$, respectively. We may assume that both functions are periodic,
otherwise we could consider $\psi_i(x) = \varphi_i(2x)$ for $i=1,2$. We can
also assume that $|a_1|^{1/b_1} < |a_2|^{1/b_2}$, and therefore
\begin{equation}
 \label{e1}
\varphi_1(x) = \varphi_2(x)c^x, \ \textrm{where}\ c>1.
\end{equation}
Finally, as $\varphi_2$ is
not identically zero, we can also suppose (by applying an appropriate
translation if needed) that $\varphi_2(0)\neq 0$.


Suppose that $b_1/b_2 \in\QQ$. Then $\varphi_1$ and $\varphi_2$
are periodic mod a common value $p$. But this is
impossible, since $c^x\neq 1$ when $x\neq 0$.

Therefore $b_1/b_2 \notin \QQ$. Then for every
(nondegenerate) interval $I\subset\RR$ there exist integers $n,k\in\ZZ$ with
$k$ arbitrarily large such that $nb_1+kb_2 \in I$. By substituting $kb_2$
into (\ref{e1}) we get $\varphi_1(kb_2) = \varphi_2(kb_2)c^{kb_2}$
for every $k\in\ZZ$, thus $\varphi_1(kb_2) = \varphi_2(0)c^{kb_2}$
for every $k\in\ZZ$. Therefore $\varphi_1(nb_1+kb_2) = \varphi_2(0)c^{kb_2}$
for every $n,k\in\ZZ$, which yields that $\varphi_1$ is unbounded on
$I$. As $I$ was arbitrary, $\varphi_1$ is unbounded on every subinterval of
$\RR$. But $\varphi_1$ is of 
Baire class~1, so it has a point of continuity (see e.g.~\cite[24.15]{Kec}), hence
it must be bounded on some interval, a contradiction.
\ep

\br
The impossibility of (\ref{e1})
is closely related to the well known statement that the
identity function is not the sum of two measurable periodic functions
(though it is surprisingly the sum of two periodic functions; see
e.g.~\cite{LR} and \cite{Kel2}). Indeed, taking the logarithm of (\ref{e1}),
we would obtain a representation of the
identity function as the sum of two Baire class 1 periodic functions;
the only problem is that our
functions can vanish at certain points.
\er

\bt\lab{diffsde}
Let $D_i f = g_i$ ($i\in I$) be a system of difference equations, and
suppose that every difference operator consists of at most two terms;
that is for every $i\in I$ the $i^{th}$ equation is of the form
\[
a_i^{(1)} f\left(x+b_i^{(1)}\right) + a_i^{(2)}
f\left(x+b_i^{(2)}\right) = g_i(x).
\]
Then if every countable subsystem has a Baire class~1
solution, then the whole system has one as well.
\et

\bp
If any of the equations consists of a single term, then it has a unique
solution, so we are clearly done. Thus, by applying a translation and 
multiplying by a real number, we may assume that every equation is of
the form
\[
f(x+b_i)-a_i f(x) = g_i(x).
\]

First suppose that $|a_{i_1}|^{1/b_{i_1}}
\neq |a_{i_2}|^{1/b_{i_2}}$ for some $i_1,i_2\in I$. Then 
it easily follows from Lemma \ref{l:ab}
that the two corresponding equations have a unique common Baire class~1
solution. This clearly solves the whole system, as every triple of equations
is solvable.

So we can assume that there exists a $c>0$ such that $|a_i|^{1/{b_i}} = c$
for every $i\in I$. If we divide the $i^{th}$ equation by $c^{x+b_i}$ and
introduce the new unknown function $\tilde{f}(x) = f(x)/c^x$, and new right
hand side $\tilde{g_i}(x) = g(x)/c^{x+b_i}$, then our equations will attain
the form (dropping the tildes) 
$\Delta_{b_i} f (x) = f(x+b_i) -  f(x) = g_i \ (i\in I^-)$ or
$f(x+b_i) +  f(x) = g_i \  
(i\in I^+)$, where $I=I^- \cup I^+$. We put 
$B^- =\{ b_i :i\in I^- \}$ and $B^+ =\{ b_i :i\in I^+ \}.$
There are countable subsets $J^- \subset I^-$ and $J^+ \subset I^+$
such that
$E^- =\{ b_i :i\in J^- \}$ is relatively dense in $B^-$, and
$E^+ =\{ b_i :i\in J^+ \} $ is relatively dense in $B^+$.


By assumption, there exists a common Baire class~1 solution $f$ to the
equations with indices $J^- \cup J^+ .$
We claim that $f$
is a solution to the whole system. 

First let $i\in I^-$. As $J^- \cup \{i\}$ is also
countable, we can choose a Baire class~1 function $f^-$ such that
$\Delta_{b_{j}} f^- = g_{j} $ for every $j\in J^-$ and $\Delta_{b_i}
f^-  = g_i .$
 
Put $f' = f-f^-$. Then for every $j\in J^-$ we have
\[
\Delta_{b_{j}} f' = \Delta_{b_{j}} (f-f^-) =
\Delta_{b_{j}}f - \Delta_{b_{j}}f^- = g_{j}-g_{j}=0,
\]
thus $f'$ is periodic mod $b_{j}$ for each $j\in J^- .$
Let $G^- =\langle E^- \rangle ;$ then $f'$ is periodic mod each
element of $G^- .$

We distinguish between two cases. If $G^-$ is dense in $\RR ,$
then $f'$ must be a constant function $c$, for
otherwise it would attain two distinct values on dense sets, so it would
have no point of continuity, which is impossible as $f'$ is Baire class~1. 

Thus
\[
\Delta_{b_i}f = \Delta_{b_i}(f^- + c) = \Delta_{b_i} f^- +\Delta_{b_i}
c = g_i + 0  = g_i,
\]
which completes the proof in the first case.

If, on the other hand, $G^-$ is not dense in $\RR$ then
$G^- =\ZZ d$ for some $d\in \RR .$ In particular, $G^-$ is discrete.
Then so is $E^-$ and thus $E^- =B^-$ as $E^-$ is dense in $B^-$. Since
$b_i\in B^- =E^- ,$ 
there is a $j\in J^-$ with $b_i =b_j $ which obviously implies $g_i
=g_j .$ Therefore, $f$ satisfies $\Delta _{b_i} f=\Delta _{b_j} f=g_j
=g_i .$

Let now $i\in I^+$. Choose a Baire class~1 function $f^+$ such that
$f^+(x+b_{j}) + f^+(x) = g_{j}(x)$ for every $j\in J^+$ and
$f^+ (x+b_i) + f(x)  = g_i(x)$ for every $x\in \RR .$
Put $f' = f-f^+$. Then $f'$ is easily seen to be anti-periodic mod
$b_{j}$, hence periodic mod  
$2b_{j}$ for every $j\in J^+$, hence it is also
periodic mod $G^+=\langle\{2b_{j} : j\in J^+ \}\rangle$.
If $G^+$ is dense in $\RR$, then $f'$ must be a constant function
$c$. But $f'$ is anti-periodic, so $c=0$. Therefore $f=f^+$, so $f$
clearly solves the $i^{th}$ equation.

On the other hand, if $G^+$ is discrete then so is $E^+$ and then we
can complete the proof as in the previous case. \ep

\bq
Is it true that $\solc(\{f : f \textrm{ is Baire class}~1 \}) = \omega_1$?
\eq

\section{Lebesgue and Baire measurable functions}\lab{s:LB}

As in Section \ref{s:s}, $\iI$ shall denote
the ideals $\iN$ or $\iM .$ Thus $\iB _{\iI}$ equals the $\sigma$-algebra
of Lebesgue or Baire measurable sets.

The goal of this section is to prove upper and lower estimates for 
$\solc(\{f: f\textrm{ is $\iB_{\iI}$-measurable}\})$ in terms of
some cardinal invariants of the ideal $\iI .$ These estimates give 
the exact value of the solvability cardinal consistently.

\bd
\begin{align*}
\add (\iI) &= \min \{|\iA| : \iA
\subset \iI, \ \bigcup \iA \notin \iI\},\\
\non (\iI) &= \min \{|A| : A\subset\RR, \ A\notin\iI\},\\
\cof (\iI) &= \min \{|\iA| : \iA \subset \iI, \ \forall I\in\iI \ \exists
A\in\iA ,\ I\subset A\}.
\end{align*}
\ed

\br\lab{r:inv}
Note that
$\omega_1 \leq \add(\iI) \le \non(\iI) \le \cof(\iI) \leq 2^\omega$.
The last
inequality follows from $\cof (\iI) \leq |\iA|$, where $ \iA= \{B\in
\iI: B \textrm{ is Borel}\}$. It is also easy to see that
$\add(\iI) \le \cf(\non(\iI))$.
\er

Before we prove our estimates (Theorems \ref{t:cof} and \ref{t:non})
we need some preparation.

\bd
For a set $H\subset\RR$ define
\[
\iD_H = \{D\in\iD : D = \sum_{i=1}^n a_i T_{b_i}, \ b_i \in H \textrm{ for
every } i=1,\ldots,n \}.
\]
\ed

\bt\lab{H}
Let $S$ be a solvable system of difference equations and
$H\subset\RR$. Then $S$ has a solution that is identically zero on $H$
if and only if whenever $(D,g)$ is deducible from $S$ and
$D\in\iD_H$ then $g(0) = 0$.
\et

\bp
The proof is again a variation of the proof of
\cite[Thm.~2.2]{La2}.

First suppose that $f$ is a solution to $S$ vanishing on $H$, and let
$(D,g)$ be deducible from $S$ such that $D\in\iD_H$; that is, $D =
\sum_{i=1}^n a_i T_{b_i}$, where $b_i \in H$ for every $i$. Then $g(0)
= (Df)(0)
= \sum_{i=1}^n a_i f(b_i) = 0$, since $b_i\in H$ for every $i$.

Suppose now that whenever $(D,g)$ is deducible from $S$ and
$D\in\iD_H$ then $g(0) = 0$. Let
\[
\iA = \{D\in\iD : \exists g \ (D,g) \textrm{ is deducible from } S
\}.
\]
Then $\iA$ is a linear subspace of $\iD$. Define
\[
L(D) = g(0) \ \ \ (D\in\iA),
\]
where $(D,g)$ is deducible from $S$. To see that $L$ is well defined
note that $S$ is solvable, and $Df=g$ whenever $(D,g)$ is deducible
from $S$ and $f$ is a solution to $S.$
$L$ is clearly linear, and by assumption
vanishes on $\iA\cap\iD_H$. Define the linear space
\[
\iB = \{A+D : A\in\iA, D\in\iD_H\}
\] 
and the linear map
\[
L'(A+D) = L(A) \ \ \ (A\in\iA, \ D\in\iD_H),
\]
which is clearly well defined and extends $L$. Moreover, $L'$
vanishes on $\iD_H$. Let $L^* : \iD\to\RR$ be a linear extension of
$L'$, and set $f(x) = L^*(T_x) \ (x\in\RR)$.  We claim that $f$ is
a solution to $S$ vanishing on $H$, which will complete the proof.

First, $L^*(T_x) = f(x) = (T_x f)(0)$ for every $x\in\RR$, and $L^*$
is linear, so $L^*(D) = (D f)(0)$ for every $D\in\iD$. If $(D,g)$ is
deducible from $S$ then so is $(T_x D,T_x g)$ for every $x\in\RR$,
hence $(Df)(x) = (T_xDf)(0) = L^*(T_xD) = (T_x g)(0) = g(x)$, so $f$
solves $S$.

Finally, $f$ vanishes on $H$, for if $x\in H$ then $T_x\in\iD_H$, so
$f(x) = L^*(T_x) = 0$, since $L^*$ vanishes on $\iD_H$.
\ep

\bt\lab{zero}
Let $S$ be a system of difference equations such that for every
$(D,g) \in S$ we have $g = 0 \ \iI$-a.e. If there exists a
$\iB_{\iI}$-measurable solution to $S$, then there is also one which is zero
$\iI$-a.e.
\et

\bp
Let $f$ be a $\iB_{\iI}$-measurable solution to $S$. 
If $\iI
= \iN$ then let $H$ be the set of points of approximate continuity of
$f$, while if $\iI = \iM$ then let $H$ be a residual (= comeager) set
on which $f$ is (relatively) continuous. Then $\RR \setminus H\in 
\iI .$ It is sufficient to show that
there exists a solution to $S$ vanishing on $H$. Using the previous
theorem we need to show that if $(D,g)$ is deducible from $S$ and
$D\in\iD_H$ then $g(0) = 0$. Let $D = \sum_{i=1}^n a_i T_{b_i}$, where
$b_i \in H$ for every $i$. As $(D,g)$ is deducible from $S$, $g=0 \ \iI$-a.e.~and
$f$ solves $(D,g)$. So $\sum_{i=1}^n a_i f(x+b_i) = g(x)$ for every
$x\in\RR$.

If $\iI = \iN$ then, using $b_i \in H$, we obtain that $g$ is
approximately continuous at $0$.  If $\iI = \iM$ then, using $b_i \in
H$, we obtain that $g$ is (relatively) continuous on the residual set
$\bigcap_{i=1}^n (H-b_i)$, which contains $0$. But in both cases $g = 0 \
\iI$-a.e., so we obtain $g(0) = 0$ as required.
\ep

\bt \label{t:cof}
$\solc(\{f: f\textrm{ is $\iB_{\iI}$-measurable}\}) \leq [\cof (\iI)]^+$.
\et

\bp
Let $S$ be such that each subsystem of cardinality at most $\cof (\iI)$
has a $\iB_{\iI}$-measurable solution. We have to show that $S$ has a
$\iB_{\iI}$-measurable solution. By Theorem \ref{ae} there exists a
$\iB_{\iI}$-measurable $f_0$ that is an 
$\iI$-a.e.~solution to $S$. Define a new system as follows.
\[
S' = \{(D,g-D f_0) : (D,g) \in S\}.
\]
Then $S'$ has a $\iB_{\iI}$-measurable solution if and only if $S$ has one,
and every subsystem of $S'$ of cardinality at most $\cof (\iI)$
has a $\iB_{\iI}$-measurable solution. Moreover, each right hand side
$g-D f_0$ is $0\ \iI$-a.e. Let
\[
S^* = \{(D,g) : (D,g) \textrm{ is deducible from } S'\}.
\]
Then $f$ solves $S'$ if and only if it solves $S^*$, moreover, each right
hand side of $S^*$ is $0\ \iI$-a.e. Also, 
every subsystem of $S^*$ of cardinality at most $\cof (\iI)$
has a $\iB_{\iI}$-measurable solution
since every $(D,g)\in S^*$ is deducible from a finite subsystem of $S' .$
In addition, every equation deducible
from $S^*$
is already in $S^*$. 

Now we prove that $S^*$ has a $\iB_{\iI}$-measurable
solution, which will complete the proof. By Theorem \ref{zero} every
subsystem of $S^*$ of cardinality at most $\cof (\iI)$
has an $\iI$-a.e.~zero solution. We claim that $S^*$ itself has such a
solution. Suppose on the contrary that this is not true. Let
$\iA\subset\iI$ be such that $|\iA| = \cof (\iI)$ and $\forall I\in\iI \ \exists
A\in\iA ,\ I\subset A$. For any $A\in\iA$ the system $S^*$ has no
solution vanishing outside $A$. By Theorem \ref{H} this means that
there exists a $(D_A,g_A) \in S^*$ such that $D_A \in
\iD_{\RR\setminus A}$ and $g_A(0) \neq 0$. 

The system $\{(D_A,g_A) : A\in\iA\}$ is of cardinality $\cof (\iI)$,
hence it has a solution $f$ vanishing $\iI$-a.e. Let $A_0 \in\iA$ be such
that $f$ vanishes outside $A_0$. Then $D_{A_0} f = g_{A_0}$, thus $(D_{A_0} f)(0) =
g_{A_0}(0) \not= 0$, but on the other hand $D_{A_0} \in \iD_{\RR \setminus
A_0}$, so $(D_{A_0} f)(0) = 0$. This contradiction finishes the proof.
\ep

\bt\lab{t:non} $\solc(\{f: f\textrm{ is $\iB_{\iI}$-measurable}\}) \geq [\cf (\non
(\iI))]^+
\ge \left[ \add (\iI ) \right] ^+ \ge \omega _2 .$
\et

\bp
We have to construct an $S$ with no $\iB_{\iI}$-measurable solutions such
that each subsystem of cardinality less than $\cf (\non (\iI))$ has a
$\iB_{\iI}$-measurable solution. 

First we construct a set $B\su\RR$ such that (i) $B\notin \iI$, $\RR\setminus B
\notin\iI$, (ii) $|(B+b)\Delta B| < \non(\iI)$ for every $b\in B$, and (iii) $B
\cap (-B) =\emptyset$.

Let $V\subset\RR$ be such that
$V\notin\iI$ and $|V| = \non (\iI)$. We may assume that $V$ is a linear
space over the rationals. Let $\{v_\al : \al<\non (\iI)\}$ be a basis of
$V$. Represent the nonzero elements of $V$ as $v = \sum_{i=1}^n q_i v_{\al_i}$,
where $q_i \in \QQ\setminus\{0\}$ and $\al_1 < \ldots <\al_n$. Define
$\varphi(v) = q_n$, and
\[
B = \{v\in V : \varphi(v) > 0\}.
\]
Clearly, (iii) holds. Note that $V= B \cup (-B) \cup \{0\}$, hence (i) is
satisfied. Let $b\in B\setminus \{0\}$ be arbitrary. Suppose
$b = \sum_{i=1}^n q_i v_{\al_i}$, where $q_i \in \QQ\setminus\{0\}$
and $\al_1 < \ldots <\al_n$. Then $(B+b) \Delta B$ is included in the linear
space generated by $\{v_\al : 
\al \leq \al_n\}$, 
which is of cardinality less than $\non (\iI)$. So (ii) holds as well. 

We claim that the system
\[
S = \{(\Delta_b,\Delta_b \chi_B) : b\in B\}
\]
satisfies the requirements.
First we check that each right hand side is zero $\iI$-a.e. 
Indeed, if $b\in B$ then $\{x\in\RR : (\Delta_b \chi_B)(x)
\neq 0\} \su (B+b) \Delta B \in\iI$, since $|(B+b) \Delta B|<\non(\iI)$. 


Suppose that $S$ has a $\iB_{\iI}$-measurable solution. Then, by Theorem
\ref{zero}, $S$ has an $\iI$-a.e.~zero solution $f_0$ as well. Then
$\Delta_b f_0 = \Delta_b \chi_B$ for every $b\in B$,
so $f_0 - \chi_B$ is periodic mod every $b\in B$. Then it is also periodic mod
each $b\in -B$. 
In particular, $f_0 - \chi_B$ is
constant on $B\cup (-B)$. 
But $f_0 = 0 \ \iI$-a.e., $B\notin\iI$, and $B\cap(-B)=\emptyset$ 
which is impossible.

What remains to show is that each subsystem $S'$ of $S$ of cardinality less
than $\cf (\non (\iI))$ has a $\iB_{\iI}$-measurable solution. Let $B'$ be the
corresponding subset of $B$, where $|B'| < \cf (\non (\iI))$. Now we put
\[
A = \langle B'\rangle +  \bigcup_{b'\in B'} \left[(B+b')\Delta B\right].
\]
Then $|A|< \non (\iI)$, hence $A \in \iI$. It is easy to see, by checking
the cases $x\in A$ and $x \notin A$, that 
$f=\chi_{B\cap A}$
is a $\iB_{\iI}$-measurable solution to $S'$.
\ep

\bcor
The Continuum Hypothesis implies
\[
\solc(\{f: f\textrm{ is measurable}\}) = \solc(\{f: f\textrm{ has the
Baire property}\}) = \omega_2 = (2^\omega)^+. 
\]
\ecor


\bq
Is $\solc(\{f: f\textrm{ is $\iB_{\iI}$-measurable}\})$ equal to $[\cof (\iI)]^+$?
Is $\solc(\{f: f\textrm{ is $\iB_{\iI}$-measurable}\})$ equal to $[\cf (\non (\iI))]^+$?
\eq



\bigskip

\noindent
\textsc{M\'arton Elekes} 

\noindent
\textsc{R\'enyi Alfr\'ed Institute} 

\noindent
\textsc{P.O. Box 127, H-1364 Budapest, Hungary}

\noindent
\textit{Email address}: \verb+emarci@renyi.hu+

\noindent
\textit{URL:} \verb+http://www.renyi.hu/~emarci+

\bigskip

\noindent
\textsc{Mikl\'os Laczkovich} 

\noindent
\textsc{E\"otv\"os Lor\'and University, Department of Analysis}
 
\noindent
\textsc{P\'az\-m\'any P\'e\-ter s\'et\'any 1/c, H-1117, Budapest, Hungary}

\noindent
\textit{Email address}: \verb+laczk@cs.elte.hu+

\bigskip

\end{document}